\documentclass[11pt]{article}
\usepackage{amsmath,amssymb,amsthm}
\usepackage{color}
\usepackage{url}
\usepackage[hidelinks]{hyperref}
 
\newtheorem{dfn}{Definition}[section]
\newtheorem{theorem}[dfn]{Theorem}
\newtheorem{lemma}[dfn]{Lemma}

\newtheorem{proposition}[dfn]{Proposition}

\def\ex{\mathrm{ex}}
\newcommand{\cI}{\mathcal I}
\newcommand{\cH}{\mathcal H}
\newcommand{\cC}{\mathcal C}
\newcommand{\cG}{\mathcal G}
\newcommand{\cL}{\mathcal L}
\newcommand{\cK}{\mathcal K}
\newcommand{\cP}{\mathcal P}
\newcommand{\one}{\mathbf 1}

\let\svthefootnote\thefootnote
\newcommand\blankfootnote[1]{%
 \let\thefootnote\relax\footnotetext{#1}%
 \let\thefootnote\svthefootnote%
}

\title{The number of $C_4$-free graphs}
\author{Tianchi Yang}
\date{}

\begin{document}
\begin{NoHyper}
\blankfootnote{School of Mathematics, Georgia Institute of Technology, Atlanta, GA 30332, US. Email: \texttt{tyang439@gatech.edu}.}
\end{NoHyper}
\maketitle

\begin{abstract}
Erd\H{o}s, Frankl and R\"odl showed that the number of graphs on $n$ vertices
with no copy of a non-bipartite graph $H$ is $2^{(1+o(1))\ex(n,H)}$, whereas
Morris and Saxton showed that this fails for $C_6$. Answering a question of
Erd\H{o}s, we prove that the number of $C_4$-free graphs on $n$ vertices is
$2^{(1+o(1))\ex(n,C_4)}$. 
The proof combines a modified hypergraph container method that treats pairs of vertices with large codegree separately, with a new bound on the number of matchings in a graph.
\end{abstract}

\section{Introduction}\label{sec:intro}

Let $H$ be a fixed graph with at least one edge. A graph is
\emph{$H$-free} if it contains no (not necessarily induced) copy of $H$,
and $\ex(n,H)$ denotes the maximum number of edges in an $H$-free graph
on $n$ vertices. How many $H$-free graphs are there on the vertex set
$[n]=\{1,\ldots,n\}$? Every subgraph of an $H$-free graph is $H$-free,
so there are at least $2^{\ex(n,H)}$ of them. On the other hand, each
of them has at most $\ex(n,H)$ edges, so their number is at most
$\sum_{m\le\ex(n,H)}\tbinom{\binom n2}{m}=2^{O(\ex(n,H)\log n)}$
whenever $\ex(n,H)=O(n^{2-\varepsilon})$ for some $\varepsilon>0$.

For non-bipartite $H$ the lower bound is sharp in the exponent.
Erd\H{o}s, Kleitman and Rothschild~\cite{EKR} proved that the number of
$K_r$-free graphs on $[n]$ is $2^{(1+o(1))\ex(n,K_r)}$, and Erd\H{o}s,
Frankl and R\"odl~\cite{EFR} extended this to every non-bipartite $H$.
Their argument uses the Szemer\'edi regularity lemma and gives no
information when $\ex(n,H)=o(n^2)$, that is, when $H$ is bipartite.

For bipartite $H$ containing a cycle, two questions arise: whether the
number of $H$-free graphs on $[n]$ is $2^{O(\ex(n,H))}$, and whether it
is even $2^{(1+o(1))\ex(n,H)}$. The first question has been answered for
many graphs. Kleitman and Winston~\cite{KW} proved that there are at most
$2^{(1.0819\ldots+o(1))n^{3/2}}$ $C_4$-free graphs on $[n]$.
Balogh and Samotij~\cite{BS11} proved the bound $2^{O(n^{2-1/s})}$ for
$K_{s,t}$-free graphs with $2\le s\le t$. The hypergraph container method,
developed independently by Balogh, Morris and Samotij~\cite{BMS} and by
Saxton and Thomason~\cite{ST}, led to much further progress. Combining it
with a balanced supersaturation theorem, Morris and Saxton~\cite{MS}
proved that there are $2^{O(n^{1+1/\ell})}$ $C_{2\ell}$-free graphs on
$[n]$ for every fixed $\ell\ge2$, which matches the Bondy--Simonovits
bound $\ex(n,C_{2\ell})=O(n^{1+1/\ell})$~\cite{BondySimonovits}.
Related results under additional assumptions can be found in~\cite{FMS,JL}.

The second question is more subtle. Morris and Saxton~\cite{MS} showed
that for some $c>0$ and infinitely many $n$ there are at least
$2^{(1+c)\ex(n,C_6)}$ $C_6$-free graphs on $[n]$. Thus the estimate
$2^{(1+o(1))\ex(n,H)}$, which holds for every non-bipartite $H$, fails for
some bipartite graphs containing a cycle.
For $C_4$, classical results give
$\ex(n,C_4)=(\frac12+o(1))n^{3/2}$~\cite{Brown,ERS,KST,Reiman};
see~\cite{MY} for sharper estimates. Erd\H{o}s asked whether the
number of $C_4$-free graphs on $[n]$ is
$2^{(1+o(1))\ex(n,C_4)}$; see Problem~3.56 in~\cite{Various99}. The
Kleitman--Winston bound leaves a constant-factor gap in the exponent,
and Balogh and Samotij~\cite{BS10} showed that almost all $C_4$-free
graphs on $[n]$ have at most $(1-\varepsilon)\ex(n,C_4)$ edges for some
fixed $\varepsilon>0$. We answer Erd\H{o}s's question affirmatively.

\begin{theorem}\label{thm:main}
The number of $C_4$-free graphs on $[n]$ is
\[
  2^{(1/2+o(1))n^{3/2}}
  =2^{(1+o(1))\ex(n,C_4)}.
\]
\end{theorem}

Our proof combines a modified hypergraph container method with a new bound
on the number of matchings in a graph. The modification is to apply the
container method only to four-cycles whose two diagonals both have small
codegree. The remaining four-cycles are deferred to the counting step inside
each container. The idea of deferring   diagonals is inspired by recent
work of Huang, Ma and the author on hypergraphs without generalized
$4$-cycles~\cite{HMY}.

Section~\ref{sec:prelim} introduces the notation, reduces
Theorem~\ref{thm:main} to sparse graphs and outlines the proof.
Section~\ref{sec:cover} constructs the containers,
Section~\ref{sec:counting} counts the $C_4$-free subgraphs of each
container, and Section~\ref{sec:matching} proves the matching bound.

\section{Preliminaries}\label{sec:prelim}

\subsection{Notation}\label{sec:notation}

All graphs are finite, simple and undirected. Unless stated otherwise,
they have vertex set $[n]$, and subgraphs of a fixed graph are spanning.
We write $N_G(v)$, $d_G(v)$ and $\Delta(G)$ for the neighborhood, degree
and maximum degree, and $e(G)=|E(G)|$. We identify graphs with their edge
sets when taking unions, intersections and differences. For distinct
vertices $a,b$, the \emph{codegree} of $a$ and $b$ is
$d_G(a,b)=|N_G(a)\cap N_G(b)|$. We write $M(G)$ for the number of
matchings in $G$, including the empty matching, and
\[
 \cI(G)=\{I\subseteq G:I\text{ is }C_4\text{-free}\}.
\]

For a hypergraph $\cH$, we write $e(\cH)$ or $|\cH|$ for its number of
hyperedges and $\cH[U]$ for its subhypergraph induced by
$U\subseteq V(\cH)$. We write $\log_2$ for base-two logarithms and $\ln$
for natural logarithms. Asymptotic notation refers to $n\to\infty$, and all
asymptotic estimates are uniform in $G$. Whenever a hypothesis
$e(G)\le cn^{3/2}$ appears, $c>0$ is fixed.
Throughout, $n$ is sufficiently large and
\begin{equation}\label{eq:parameters}
 t=\bigl\lceil\sqrt{\log_2 n}\,\bigr\rceil,\qquad
 r=t^{16},\qquad \theta=t^{-5}.
\end{equation}

Let $\cH(G)$ be the $4$-uniform hypergraph on $E(G)$ whose hyperedges
are the edge sets of the four-cycles in $G$. Then $I\in\cI(G)$ if and only
if $E(I)$ is independent in $\cH(G)$. In particular, a hypergraph
container is a set of edges of $G$, which we identify with a spanning
subgraph of $G$.

The \emph{diagonals} of a four-cycle $abcda$ are the pairs $ac$ and
$bd$. To isolate pairs with large common neighborhoods, for
$\theta>0$ define the \emph{diagonal graph} $D_\theta(G)$ on $[n]$ by
\[
 ab\in E(D_\theta(G))
 \quad\Longleftrightarrow\quad d_G(a,b)\ge\theta\sqrt n.
\]
Its edges need not be edges of $G$. For a graph $D$ on $[n]$ and an
integer $t\ge2$, let $D^t$ be the graph on $[n]$ in which two distinct
vertices are adjacent if their distance in $D$ is at most $t$.

Given a graph $D$ on $[n]$, let $\cH_D(G)\subseteq\cH(G)$ consist of
the four-cycles of $G$ neither of whose diagonals is an edge of $D$.
Then, for every $I\subseteq G$,
\begin{equation}\label{eq:constraint-decomposition}
 I\in\cI(G)
 \quad\Longleftrightarrow\quad
 \begin{cases}
 E(I)\text{ is independent in }\cH_D(G),\\
 d_I(a,b)\le1\text{ for every }ab\in E(D).
 \end{cases}
\end{equation}
We apply the container method to $\cH_D(G)$, while retaining the second
condition in~\eqref{eq:constraint-decomposition} for the subsequent
counting argument.

\subsection{Reduction to sparse graphs}\label{sec:sparse-reduction}

We use Theorem~1.2 of Morris and Saxton~\cite{MS} in the case $\ell=2$.

\begin{theorem}[\cite{MS}]\label{thm:sparse-cover}
For every $\ell\ge2$ and $\delta>0$ there is a constant $c=c(\delta,\ell)$
such that the following holds for all sufficiently large $n$. There is a
collection $\cG$ of at most $2^{\delta n^{1+1/\ell}}$ graphs on $[n]$ such
that $e(G)\le cn^{1+1/\ell}$ for every $G\in\cG$, and every
$C_{2\ell}$-free graph on $[n]$ is a subgraph of some $G\in\cG$.
\end{theorem}

The main work of the paper is the proof of the following result.

\begin{theorem}\label{thm:sparse-graph}
For every fixed $c>0$, every graph $G$ on $[n]$ with
$e(G)\le cn^{3/2}$ satisfies
\begin{equation}\label{eq:sparse-graph}
 \log_2|\cI(G)|\le\Bigl(\frac12+o(1)\Bigr)n^{3/2}.
\end{equation}
\end{theorem}

\begin{proof}[Proof of Theorem~\ref{thm:main} assuming Theorem~\ref{thm:sparse-graph}]
The lower bound follows by taking all subgraphs of an extremal $C_4$-free
graph. For the upper bound, fix $\delta>0$ and let $\cG$ be the
collection given by Theorem~\ref{thm:sparse-cover} with $\ell=2$.
Applying Theorem~\ref{thm:sparse-graph} with $c=c(\delta,2)$ gives
\[
 |\cI(K_n)|
 \le\sum_{G\in\cG}|\cI(G)|
 \le2^{\delta n^{3/2}}\cdot2^{(1/2+o(1))n^{3/2}}
 =2^{(1/2+\delta+o(1))n^{3/2}}.
\]
Since $\delta>0$ is arbitrary, the upper bound follows.
\end{proof}

Conversely, Theorem~\ref{thm:main} implies Theorem~\ref{thm:sparse-graph}
because $\cI(G)\subseteq\cI(K_n)$, so the two statements are equivalent.

\subsection{Outline of the proof}\label{sec:strategy}

Fix a graph $G$ with $e(G)\le cn^{3/2}$.  
Lemmas~\ref{lem:degree-reduction} and~\ref{lem:diagonal-reduction} give
a subgraph $G_*\subseteq G$ such that $\Delta(G_*)\le t\sqrt n$, the
graph $D=D_\theta(G_*)$ satisfies $\Delta(D)<r$, and
\[
 |\cI(G)|\le2^{o(n^{3/2})}|\cI(G_*)|.
\]
Keeping $D$ fixed, we apply the container theorem to $\cH_D(G_*)$, whose
hyperedges are the four-cycles of $G_*$ with neither diagonal in $D$.
Lemma~\ref{lem:selective-containers} gives a collection $\cC$ of at most
$2^{o(n^{3/2})}$ subgraphs of $G_*$ covering $\cI(G_*)$, with
$|\cH_D(C)|\le n^2/t$ for every $C\in\cC$. Consequently,
\[
 |\cI(G_*)|\le2^{o(n^{3/2})}\max_{C\in\cC}|\cI(C)|.
\]

A container may still contain many four-cycles with a diagonal in $D$,
so we cannot simply take all of its subgraphs. Instead, for each
$C\in\cC$ we partition every neighborhood $N_C(u)$ into $D$-connected
sets of at most $t$ vertices and split each vertex $u$ accordingly,
obtaining an auxiliary graph $H$ (Section~\ref{sec:partitions}). By
Lemma~\ref{lem:matching-representation}, each $I\in\cI(C)$ is
determined by a record of common neighbors of the pairs adjacent in
$D^t$, for which there are $2^{o(n^{3/2})}$ choices, together with a
matching in $H$. The number of matchings is bounded by the following
inequality, which we prove in Section~\ref{sec:matching}.

\begin{proposition}\label{prop:matching-degree}
For every graph $F$,
\[
 \log_2M(F)\le\frac12\sum_{z\in V(F)}\sqrt{d_F(z)}.
\]
\end{proposition}

In $H$, the copy of $u$ corresponding to a part $S$ has degree $|S|$.
Lemma~\ref{lem:component-sum}, a weighted form of the Cauchy--Schwarz
argument for $\ex(n,C_4)$, shows that the sum of $\sqrt{|S|}$ over all
parts is at most $(1+o(1))n^{3/2}$. Hence
\[
 |\cI(C)|\le2^{o(n^{3/2})}M(H)\le2^{(1/2+o(1))n^{3/2}}.
\]
The coefficient $\frac12$ in Proposition~\ref{prop:matching-degree} is
sharp already for a single edge, and it is precisely the constant that
appears in Theorem~\ref{thm:main}. Combining these bounds proves
Theorem~\ref{thm:sparse-graph}.

\section{Constructing the containers}\label{sec:cover}

In this section we prove the following proposition. Throughout, $c>0$
is fixed, $n$ is sufficiently large, and the parameters are as
in~\eqref{eq:parameters}.

\begin{proposition}\label{prop:container-preparation}
Let $G$ be a graph with $e(G)\le cn^{3/2}$. There exist a graph $D$ with
$\Delta(D)<r$ and a nonempty collection $\cC$ of at most $2^{o(n^{3/2})}$
subgraphs of $G$ such that $|\cH_D(C)|\le n^2/t$ for every $C\in\cC$ and
\begin{equation}\label{eq:prepared-count}
 |\cI(G)|\le2^{o(n^{3/2})}\max_{C\in\cC}|\cI(C)|.
\end{equation}
\end{proposition}

The proof has three steps. We first delete all edges incident to vertices
of degree larger than $t\sqrt n$ (Lemma~\ref{lem:degree-reduction}), then
delete further edges until the diagonal graph has maximum degree less
than $r$ (Lemma~\ref{lem:diagonal-reduction}), and finally apply the
container theorem to the four-cycles with no diagonal in the diagonal
graph (Lemma~\ref{lem:selective-containers}). In the first two steps,
the deleted edges are covered by few vertices, and their $C_4$-free
subgraphs are counted by the following lemma.

\subsection{Bounding the maximum degree}\label{sec:degree-deletions}

\begin{lemma}\label{lem:small-vertex-cover}
There is an absolute constant $c_1>0$ with the following property. Let
$G$ be a graph on an arbitrary finite vertex set, and let
$S\subseteq V(G)$ be such that every edge of $G$ has an endpoint in $S$.
Then
\[
 \log_2|\cI(G)|\le c_1|V(G)|\sqrt{|S|}.
\]
\end{lemma}

\begin{proof}
Kleitman and Winston~\cite{KW} proved that
$\log_2|\cI(K_\ell)|\le(1.0819\ldots+o(1))\ell^{3/2}$ as $\ell\to\infty$.
Hence there is an absolute constant $c_0>0$ such that
$\log_2|\cI(K_\ell)|\le c_0\ell^{3/2}$ for every integer $\ell\ge1$. We
show that $c_1=4c_0$ works.

If $S=\varnothing$, then $G$ has no edges and the claim is trivial.
Otherwise, partition $V(G)\setminus S$ into sets $W_1,\ldots,W_k$ of
size at most $|S|$, where $k\le|V(G)|/|S|$; if $S=V(G)$, take $k=1$ and
$W_1=\varnothing$. Every edge of $G$ either lies inside $S$ or joins $S$
to some $W_i$, so it lies in $G[S\cup W_i]$ for some $i$. Hence
$I\mapsto(I\cap G[S\cup W_i])_{i\in[k]}$ is injective on $\cI(G)$, and
since $|S\cup W_i|\le2|S|$,
\[
 \log_2|\cI(G)|
 \le\sum_{i=1}^{k}\log_2|\cI(G[S\cup W_i])|
 \le c_0k(2|S|)^{3/2}
 \le 4c_0|V(G)|\sqrt{|S|}.\qedhere
\]
\end{proof}

\begin{lemma}\label{lem:degree-reduction}
Let $G$ be a graph with $e(G)\le cn^{3/2}$. There is a subgraph
$G_0\subseteq G$ with $\Delta(G_0)\le t\sqrt n$ such that
\[
 |\cI(G)|\le2^{o(n^{3/2})}|\cI(G_0)|.
\]
\end{lemma}

\begin{proof}
Let $W=\{v:d_G(v)>t\sqrt n\}$, and let $G_0$ be obtained from $G$ by
deleting all edges meeting $W$. Then $\Delta(G_0)\le t\sqrt n$, and
summing degrees gives $|W|\le2e(G)/(t\sqrt n)\le2cn/t$. Every
$I\in\cI(G)$ is the edge-disjoint union of the $C_4$-free graphs
$I\cap G_0$ and $I\cap(G\setminus G_0)$. Every edge of $G\setminus G_0$
meets $W$, so Lemma~\ref{lem:small-vertex-cover} gives
\[
 |\cI(G)|
 \le|\cI(G\setminus G_0)|\,|\cI(G_0)|
 \le2^{c_1n\sqrt{2cn/t}}\,|\cI(G_0)|
 =2^{o(n^{3/2})}|\cI(G_0)|.\qedhere
\]
\end{proof}

\subsection{Bounding the maximum degree of the diagonal graph}
\label{sec:diagonal-preparation}

\begin{lemma}\label{lem:diagonal-reduction}
Let $G_0$ be a graph with $e(G_0)\le cn^{3/2}$ and
$\Delta(G_0)\le t\sqrt n$. There is a subgraph $G_*\subseteq G_0$ such
that $\Delta(D_\theta(G_*))<r$ and
\[
 |\cI(G_0)|\le2^{o(n^{3/2})}|\cI(G_*)|.
\]
\end{lemma}

\begin{proof}
We construct graphs $G_0\supseteq G_1\supseteq\cdots$ as follows. Given
$G_i$, stop if $\Delta(D_\theta(G_i))<r$. Otherwise, choose a vertex $a$
of degree at least $r$ in $D_\theta(G_i)$ and a set
$S_i\subseteq N_{D_\theta(G_i)}(a)$ with $|S_i|=r$.
Let $Z_i$ consist of all edges of $G_i$ with one endpoint in $S_i$
and the other in $N_{G_i}(a)$, and let $G_{i+1}=G_i\setminus Z_i$.

Each $b\in S_i$ has at least $\theta\sqrt n$ neighbors in
$N_{G_i}(a)$, since $ab\in E(D_\theta(G_i))$.
Note that $S_i$ and $N_{G_i}(a)$ may overlap.
Then 
$2|Z_i|\ge \sum_{b\in S_i}\theta\sqrt n = r\theta\sqrt n .$
The sets $Z_i$ are pairwise disjoint subsets of $E(G_0)$,
so the process stops after some number $k$ of steps, where
\[
 k\le\frac{2e(G_0)}{r\theta\sqrt n}
 \le\frac{2cn}{r\theta}.
\]
Let $G_*=G_k$. Then $\Delta(D_\theta(G_*))<r$, and $E(G_0)$ is the
disjoint union of $E(G_*),Z_0,\ldots,Z_{k-1}$. Every $I\in\cI(G_0)$ is
determined by its intersections with these sets, each of which is
$C_4$-free, so
\[
 |\cI(G_0)|\le|\cI(G_*)|\prod_{i=0}^{k-1}|\cI(Z_i)|,
\]
where we regard $Z_i$ as a graph on its set of endpoints. Every edge of
$Z_i$ meets $S_i$, and $Z_i$ has at most
$|S_i|+\Delta(G_0)\le r+t\sqrt n\le2t\sqrt n$ vertices. By
Lemma~\ref{lem:small-vertex-cover},
$\log_2|\cI(Z_i)|\le2c_1t\sqrt{nr}$ for every $i$. Consequently,
\[
 \log_2\frac{|\cI(G_0)|}{|\cI(G_*)|}
 \le k\cdot2c_1t\sqrt{nr}
 \le\frac{4cc_1t}{\theta\sqrt r}\,n^{3/2}
 =\frac{4cc_1}{t^2}\,n^{3/2}
 =o(n^{3/2}).\qedhere
\]
\end{proof}

\subsection{Applying the container theorem}\label{sec:fixed-d-containers}

We use the $4$-uniform case of the container theorem of Saxton and
Thomason~\cite{ST}, in the form of Definition~4.1 and Theorem~4.2 of
Morris and Saxton~\cite{MS}. For a $4$-uniform hypergraph $\cH$ and
$\sigma\subseteq V(\cH)$, let $d_\cH(\sigma)$ be the number of hyperedges
containing $\sigma$. For $v\in V(\cH)$ and $j\in\{2,3,4\}$, let
\[
 d^{(j)}(v)=\max\bigl\{d_\cH(\sigma):
 v\in\sigma\subseteq V(\cH),\ |\sigma|=j\bigr\},
\]
with $d^{(j)}(v)=0$ if $|V(\cH)|<j$. For $e(\cH)>0$ and $\tau>0$, the
\emph{codegree function} is
\[
 \delta(\cH,\tau)
 =\frac1{e(\cH)}\sum_{j=2}^{4}\frac1{\tau^{j-1}}
   \sum_{v\in V(\cH)}d^{(j)}(v).
\]
Writing $\Delta_j(\cH)=\max_{|\sigma|=j}d_\cH(\sigma)$, we have, for a
hypergraph $\cH$ with $N$ vertices,
\begin{equation}\label{eq:codegree-maximum-bound}
 \delta(\cH,\tau)
 \le\frac{N}{e(\cH)}\left(\frac{\Delta_2(\cH)}{\tau}
   +\frac{\Delta_3(\cH)}{\tau^2}+\frac{\Delta_4(\cH)}{\tau^3}\right).
\end{equation}

\begin{theorem}[\cite{MS,ST}]\label{thm:container-step}
There is $\delta_0>0$ such that the following holds for every
$0<\delta<\delta_0$. Let $\cH$ be a $4$-uniform hypergraph on $N$
vertices with $e(\cH)>0$, and let $0<\tau<1/2$ satisfy
$\delta(\cH,\tau)\le\delta$. Then there is a collection $\cC$ of subsets
of $V(\cH)$ with
\[
 |\cC|\le\exp\left(\frac{\tau\ln(1/\tau)N}{\delta}\right)
\]
such that every independent set of $\cH$ is contained in some
$C\in\cC$, and $e(\cH[C])\le(1-\delta)e(\cH)$ for every $C\in\cC$.
\end{theorem}

We apply Theorem~\ref{thm:container-step} repeatedly with the diagonal
graph $D$ fixed.

\begin{lemma}\label{lem:selective-containers}
Let $G_*$ be a graph with $e(G_*)\le cn^{3/2}$, and let
$D=D_\theta(G_*)$. There is a nonempty collection $\cC$ of at most
$2^{o(n^{3/2})}$ subgraphs of $G_*$ such that every
$I\in\cI(G_*)$ is a subgraph of some $C\in\cC$, and
$|\cH_D(C)|\le n^2/t$ for every $C\in\cC$. In particular,
\[
 |\cI(G_*)|\le2^{o(n^{3/2})}\max_{C\in\cC}|\cI(C)|.
\]
\end{lemma}

\begin{proof}
Fix $0<\delta<\delta_0$, where $\delta_0$ is given by
Theorem~\ref{thm:container-step}, and let $\tau=t^{-3}$.
We first describe one refinement step.
Let $C\subseteq G_*$. If $|\cH_D(C)|\le n^2/t$, let
$\cC(C)=\{C\}$. Otherwise, let $\cH=\cH_D(C)$, so that
$N=|V(\cH)|=e(C)\le cn^{3/2}$ and $e(\cH)>n^2/t$.

We bound the codegrees of $\cH$. Consider two edges of $C$. If they
share a vertex, say they are $ua$ and $ub$, then every four-cycle of
$C$ containing both has the form $uavbu$ and has $ab$ as a diagonal.
Hence no hyperedge of $\cH$ contains both edges if $ab\in E(D)$, and
otherwise there are fewer than $\theta\sqrt n$ choices for $v$. Two
disjoint edges lie together in at most two four-cycles, and three edges
lie together in at most one. Therefore
\[
 \Delta_2(\cH)\le\theta\sqrt n+2,\qquad
 \Delta_3(\cH),\,\Delta_4(\cH)\le1,
\]
and~\eqref{eq:codegree-maximum-bound} and~\eqref{eq:parameters} give
\[
 \delta(\cH,\tau)
 \le\frac{cn^{3/2}t}{n^2}\left(\frac{\theta\sqrt n+2}{\tau}
   +\frac1{\tau^2}+\frac1{\tau^3}\right)
 \le c\left(\frac1t+\frac{4t^{10}}{\sqrt n}\right)=o(1).
\]
Thus $\delta(\cH,\tau)\le\delta$ and $0<\tau<1/2$ for $n$ sufficiently
large. Every $I\in\cI(C)$ is independent in $\cH$, so
Theorem~\ref{thm:container-step} gives a collection $\cC(C)$ of
subgraphs of $C$ covering $\cI(C)$. Since $D$ is fixed, the hyperedges
of $\cH$ contained in $E(C')$ are precisely the members of $\cH_D(C')$.
Hence each $C'\in\cC(C)$ satisfies
\[
 |\cH_D(C')|\le(1-\delta)|\cH_D(C)|.
\]
In both cases,
$\cC(C)$ covers $\cI(C)$,
\begin{equation}\label{eq:one-step}
 |\cH_D(C')|\le\max\Bigl\{(1-\delta)|\cH_D(C)|,\frac{n^2}{t}\Bigr\}
 \quad\text{for every }C'\in\cC(C),
\end{equation}
and, since $\tau\log_2(1/\tau)=3t^{-3}\log_2t$,
\begin{equation}\label{eq:one-step-size}
 \log_2|\cC(C)|\le\frac{e(C)\tau\log_2(1/\tau)}{\delta}
 \le\frac{3c}{\delta}\cdot\frac{\log_2t}{t^3}\,n^{3/2}.
\end{equation}

Now let $\cC_0=\{G_*\}$ and
$\cC_{i+1}=\bigcup_{C\in\cC_i}\cC(C)$ for $i\ge0$. By induction on $i$,
using that every $I\in\cI(G_*)$ contained in $C$ belongs to $\cI(C)$,
the collection $\cC_i$ covers $\cI(G_*)$. Since
$|\cH_D(G_*)|\le n^4$, induction using~\eqref{eq:one-step}
and~\eqref{eq:one-step-size} also gives
\[
 |\cH_D(C)|\le\max\Bigl\{(1-\delta)^in^4,\frac{n^2}{t}\Bigr\}
 \quad\text{for }C\in\cC_i,
 \qquad
 \log_2|\cC_i|\le\frac{3ci}{\delta}\cdot\frac{\log_2t}{t^3}\,n^{3/2}.
\]
Let $m=\lceil\ln(n^2t)/(-\ln(1-\delta))\rceil$ and $\cC=\cC_m$. Then
$(1-\delta)^mn^4\le n^2/t$, so $|\cH_D(C)|\le n^2/t$ for every
$C\in\cC$. Moreover, $m=O(\log n)=O(t^2)$, so
\[
 \log_2|\cC|=O\Bigl(\frac{\log_2t}{t}\,n^{3/2}\Bigr)=o(n^{3/2}).
\]
The collection $\cC$ is nonempty because it covers the empty graph, and
\[
 |\cI(G_*)|\le\sum_{C\in\cC}|\cI(C)|
 \le2^{o(n^{3/2})}\max_{C\in\cC}|\cI(C)|.\qedhere
\]
\end{proof}

\begin{proof}[Proof of Proposition~\ref{prop:container-preparation}]
Apply Lemmas~\ref{lem:degree-reduction} and
\ref{lem:diagonal-reduction} successively to obtain graphs
$G_*\subseteq G_0\subseteq G$. Then $\Delta(D_\theta(G_*))<r$, and
Lemma~\ref{lem:selective-containers} applied to $G_*$ yields a
collection $\cC$ such that
$|\cH_{D_\theta(G_*)}(C)|\le n^2/t$ for every $C\in\cC$.
Combining the counting bounds in the three lemmas gives
\[
 |\cI(G)|
 \le2^{o(n^{3/2})}|\cI(G_0)|
 \le2^{o(n^{3/2})}|\cI(G_*)|
 \le2^{o(n^{3/2})}\max_{C\in\cC}|\cI(C)|.
\]
Thus $D_\theta(G_*)$ and $\cC$ satisfy the proposition.
\end{proof}

\section{Counting the \texorpdfstring{$C_4$}{C4}-free subgraphs of a container}
\label{sec:counting}

In this section we prove the following proposition, which together with
Proposition~\ref{prop:container-preparation} implies
Theorem~\ref{thm:sparse-graph}.

\begin{proposition}\label{prop:matching-count}
Let $c>0$ be fixed. If $G$ and $D$ are graphs on $[n]$ satisfying
\[
 e(G)\le cn^{3/2},\qquad \Delta(D)\le r,\qquad |\cH_D(G)|\le\frac{n^2}{t},
\]
we have
\[
 \log_2|\cI(G)|\le\Bigl(\frac12+o(1)\Bigr)n^{3/2}.
\]
\end{proposition}

\begin{proof}[Proof of Theorem~\ref{thm:sparse-graph}]
Let $G$ be a graph on $[n]$ with $e(G)\le cn^{3/2}$, and let $D$ and
$\cC$ be given by Proposition~\ref{prop:container-preparation}. Every
$C\in\cC$ satisfies $e(C)\le e(G)\le cn^{3/2}$ and
$|\cH_D(C)|\le n^2/t$, and $\Delta(D)<r$. By~\eqref{eq:prepared-count}
and Proposition~\ref{prop:matching-count} applied to each $C\in\cC$,
\[
 |\cI(G)|\le2^{o(n^{3/2})}\max_{C\in\cC}|\cI(C)|
 \le2^{(1/2+o(1))n^{3/2}}.\qedhere
\]
\end{proof}

In Section~\ref{sec:partitions} we encode each $C_4$-free subgraph of
$G$ by a record of common neighbors and a matching in an auxiliary
graph $H$. In Section~\ref{sec:matching-count-proof} we bound the number
of records directly, and the number of matchings by combining
Proposition~\ref{prop:matching-degree} with a weighted form of the
Cauchy--Schwarz argument for $\ex(n,C_4)$.

\subsection{A matching representation}\label{sec:partitions}

Let $G$ and $D$ be graphs on $[n]$. For each $u\in[n]$, start with
$N_G(u)$ and repeatedly remove a $t$-vertex set that induces a connected
subgraph of $D$, for as long as such a set exists. Let $\cL_u$ be the
family of sets removed in this way, and let $\cK_u$ be the family of
vertex sets of the connected components of the remaining graph.

Every member of $\cL_u$ has size exactly $t$. Moreover, every member
of $\cK_u$ has size less than $t$. Otherwise a remaining component with
at least $t$ vertices would contain a connected $t$-vertex subgraph, and
the procedure could continue. Thus
\[
 \cP_u=\cL_u\cup\cK_u
\]
is a partition of $N_G(u)$ into sets of size at most $t$, each inducing
a connected subgraph of $D$. Since the members of $\cK_u$ are distinct
components of the remaining graph,
\begin{equation}\label{eq:K-separated}
 \text{no edge of $D$ joins two distinct members of $\cK_u$.}
\end{equation}

We next split each vertex according to the partition of its neighborhood.
Let $H$ be the graph with vertex set
\[
 V(H)=\{u_S: u\in[n],\ S\in\cP_u\}.
\]
For each edge $uv\in E(G)$, let $S\in\cP_u$ be the part containing $v$, let
$T\in\cP_v$ be the part containing $u$, and set $\varphi(uv)=u_Sv_T$. The edge
set of $H$ is $\{\varphi(e):e\in E(G)\}$. Since $u_Sv_T$ determines $uv$,
the map $\varphi:E(G)\to E(H)$ is a bijection.
By construction, the edges of $H$ at $u_S$ are exactly
$\varphi(uv)$ with $v\in S$, so
\begin{equation}\label{eq:detachment-degree}
 d_H(u_S)=|S|.
\end{equation}
For $F\subseteq G$, let $\varphi(F)$ be the spanning subgraph of $H$
with edge set $\varphi(E(F))$, and define $\varphi^{-1}$ on spanning
subgraphs of $H$ accordingly.

Now let $I\in\cI(G)$. For $ab\in E(D^t)$ with $a<b$, let $Q_I(ab)$ be
the common neighbor of $a$ and $b$ in $I$ if it exists, and
$Q_I(ab)=\ast$ otherwise. Since $I$ is $C_4$-free, $a$ and $b$ have at
most one common neighbor in $I$, so $Q_I:E(D^t)\to[n]\cup\{\ast\}$ is
well defined. Let $I_1$ be the graph on $[n]$ with edge set
\[
 E(I_1)=\{aw:ab\in E(D^t),\ a<b,\ Q_I(ab)=w\ne\ast\},
\]
and let $I_2=I\setminus I_1$ and $M=\varphi(I_2)$.
\begin{lemma}\label{lem:matching-representation}
For every $I\in\cI(G)$, the graph $M=\varphi(I_2)$ is a matching in $H$. Moreover,
\[
 |\cI(G)|\le(n+1)^{e(D^t)}M(H).
\]
\end{lemma}

\begin{proof}
We first show that $N_{I_2}(u)$ is an independent set in $D^t$ for every
$u\in[n]$. Suppose that $a<b$, $ab\in E(D^t)$ and $a,b\in N_{I_2}(u)$.
Then $u$ is a common neighbor of $a$ and $b$ in $I$, so $Q_I(ab)=u$ and
$au\in E(I_1)$, contradicting $au\in E(I_2)$.

Now let $S\in\cP_u$. Since $D[S]$ is connected and $|S|\le t$, any two
vertices of $S$ are joined by a path in $D[S]$ of length at most
$|S|-1<t$, so $S$ is a clique in $D^t$. Hence $|N_{I_2}(u)\cap S|\le1$.
The edges of $M$ at $u_S$ are the images of the edges of $I_2$ from $u$
to $S$, so $d_M(u_S)=|N_{I_2}(u)\cap S|\le1$, and $M$ is a matching.

To prove the inequality, note that $I=I_1\cup\varphi^{-1}(M)$ and $I_1$ is determined by $Q_I$, so the map $I\mapsto(Q_I,M)$ is injective on $\cI(G)$. There are at most $(n+1)^{e(D^t)}$ choices for $Q_I$ and at most $M(H)$ choices for $M$, which gives \[ |\cI(G)|\le(n+1)^{e(D^t)}M(H).\qedhere \]
\end{proof}

\subsection{Proof of Proposition~\ref{prop:matching-count}}
\label{sec:matching-count-proof}\label{sec:weighted}

Let $G$ and $D$ satisfy the hypotheses of
Proposition~\ref{prop:matching-count}, and let $\cL_u$, $\cK_u$ and $H$
be as in Section~\ref{sec:partitions}. By
Lemma~\ref{lem:matching-representation},
\begin{equation}\label{eq:record-matching}
 \log_2|\cI(G)|\le e(D^t)\log_2(n+1)+\log_2M(H).
\end{equation}
We show that the first term is $o(n^{3/2})$ and the second is at most
$(\frac12+o(1))n^{3/2}$.
 
Since $\Delta(D)\le r$ and $r\ge2$, every vertex has at most
$\sum_{j=1}^tr^j\le2r^t$ vertices at distance at most $t$ in $D$, so
$e(D^t)\le nr^t$. By~\eqref{eq:parameters},
$\log_2(r^t)=16t\log_2t=o(\log n)$, so $r^t=n^{o(1)}$ and
\[
 e(D^t)\log_2(n+1)\le nr^t\log_2(n+1)=n^{1+o(1)}=o(n^{3/2}).
\]
Since the sets in $\cL_u$ are disjoint $t$-subsets of $N_G(u)$, \[ \frac12\sum_{u\in[n]}\sum_{S\in\cL_u}\sqrt{|S|} =\frac1{2\sqrt t}\sum_{u\in[n]}\sum_{S\in\cL_u}|S| \le\frac1{2\sqrt t}\sum_{u\in[n]}d_G(u) =\frac{e(G)}{\sqrt t} =o(n^{3/2}). \]
Hence Proposition~\ref{prop:matching-degree} and~\eqref{eq:detachment-degree} give
\begin{equation}\label{eq:matching-reduction} \log_2 M(H) \le  o(n^{3/2}) +\frac12\sum_{u\in[n]}\sum_{S\in\cK_u}\sqrt{|S|}. \end{equation}

To bound the second term,  we adapt the classical Cauchy--Schwarz argument for $\ex(n,C_4)$, applied to normalized indicator vectors of the sets in $\cK_u$. 

\begin{lemma}\label{lem:component-sum}
Let $G$ and $D$ be graphs on $[n]$ with $\Delta(D)\le r$. For each
$u\in[n]$, let $\cK_u$ be a family of pairwise disjoint nonempty subsets
of $N_G(u)$ such that no edge of $D$ joins two distinct members of
$\cK_u$. Then
\[
 \sum_{u\in[n]}\sum_{S\in\cK_u}\sqrt{|S|}
 \le\sqrt{n\bigl(n(n-1)+2(r+1)e(G)+4|\cH_D(G)|\bigr)}.
\]
\stepcounter{equation}
\end{lemma}

\begin{proof}
For $S\subseteq[n]$, let $\one_S\in\mathbb R^n$ be its indicator vector,
and let $\one=\one_{[n]}$. For $u\in[n]$, let
\[
 f_u=\sum_{S\in\cK_u}\frac{\one_S}{\sqrt{|S|}}.
\]
Since the members of $\cK_u$ are disjoint subsets of $N_G(u)$,
\[
 \langle\one,f_u\rangle=\sum_{S\in\cK_u}\sqrt{|S|},\qquad
 \|f_u\|_2^2=|\cK_u|\le d_G(u),\qquad
 0\le f_u\le\one_{N_G(u)}\ \text{coordinatewise}.
\]
By the Cauchy--Schwarz inequality,
\[
 \Bigl(\sum_u\sum_{S\in\cK_u}\sqrt{|S|}\Bigr)^2
 =\Bigl\langle\one,\sum_uf_u\Bigr\rangle^2
 \le n\Bigl\|\sum_uf_u\Bigr\|_2^2
 =n\Bigl(\sum_u\|f_u\|_2^2+2\sum_{\{a,b\}\in\binom{[n]}2}\langle f_a,f_b\rangle\Bigr),
\]
Here  $\sum_u\|f_u\|_2^2\le\sum_ud_G(u)=2e(G)$.  It therefore suffices to
show that
\begin{equation}\label{eq:cross-terms}
 \sum_{ab\in E(D)}\langle f_a,f_b\rangle\le r\,e(G)
 \qquad\text{and}\qquad
 \sum_{ab\notin E(D)}\langle f_a,f_b\rangle\le\binom n2+2|\cH_D(G)|.
\end{equation}

\smallskip\noindent\emph{Case 1: $ab\in E(D)$.} Since
$0\le f_a\le\one_{N_G(a)}$ and $0\le f_b\le\one_{N_G(b)}$
coordinatewise,
\[
 \langle f_a,f_b\rangle\le
 \langle\one_{N_G(a)},\one_{N_G(b)}\rangle=d_G(a,b).
\]
Counting the pairs $(ab,v)$ with $ab\in E(D)$ and $v$ a common neighbor
of $a$ and $b$ in $G$, and using $\Delta(D)\le r$, we obtain
\[
 \sum_{ab\in E(D)}\langle f_a,f_b\rangle
 \le\sum_{ab\in E(D)}d_G(a,b)
 =\sum_{v\in[n]}e\bigl(D[N_G(v)]\bigr)
 \le\sum_{v\in[n]}\frac{r\,d_G(v)}2=r\,e(G).
\]

\smallskip\noindent\emph{Case 2: $ab\notin E(D)$.} 
Let $W$ be the set of vertices covered both by
$\cK_a$ and by $\cK_b$. If $W=\varnothing$, then
$\langle f_a,f_b\rangle=0$. Otherwise, fix $x\in W$, let $A\in\cK_a$
and $B\in\cK_b$ be the members containing $x$, and let $S=A\cap B$. The
product $f_a(v)f_b(v)$ vanishes for $v\notin W$, equals
$1/\sqrt{|A||B|}$ for $v\in S$, and is at most $1$ for
$v\in W\setminus S$. Since $|S|\le\min\{|A|,|B|\}\le\sqrt{|A||B|}$,
\[
 \langle f_a,f_b\rangle
 \le\frac{|S|}{\sqrt{|A||B|}}+|W\setminus S|
 \le1+|W\setminus S|.
\]
Let $y\in W\setminus S$. Then $y\notin A$ or $y\notin B$, and by symmetry we may assume
$y\notin A$. As $y\in W$, it lies in a member of $\cK_a$ other than $A$,
so $xy\notin E(D)$ by hypothesis.  
Then $axbya$ is a four-cycle in $\cH_D(G)$ with diagonal $ab$, and distinct choices of $y$ give distinct four-cycles. Hence $\langle f_a,f_b\rangle$ is at most $1$ plus the number of four-cycles in $\cH_D(G)$ with diagonal $ab$. Every four-cycle has two diagonals, so summing over all $ab\notin E(D)$ gives the second bound in~\eqref{eq:cross-terms}.

Combining these estimates,
\[
 \Bigl\|\sum_uf_u\Bigr\|_2^2
 \le2e(G)+2r\,e(G)+n(n-1)+4|\cH_D(G)|,
\]
This proves Lemma~\ref{lem:component-sum}.
\end{proof}

\begin{proof}[Completion of the proof of Proposition~\ref{prop:matching-count}]
By~\eqref{eq:K-separated}, Lemma~\ref{lem:component-sum} applies to the
families $\cK_u$. By the hypotheses and~\eqref{eq:parameters},
\[
 (r+1)e(G)=O(t^{16}n^{3/2})=o(n^2),\qquad
 |\cH_D(G)|\le\frac{n^2}{t}=o(n^2),
\]
so~Lemma~\ref{lem:component-sum} gives
\[
 \sum_{u\in[n]}\sum_{S\in\cK_u}\sqrt{|S|}
 \le\sqrt{n\bigl(n^2+o(n^2)\bigr)}=(1+o(1))\,n^{3/2}.
\]
By~\eqref{eq:matching-reduction},
$\log_2M(H)\le(\frac12+o(1))n^{3/2}$. Together with the bound on the
records and~\eqref{eq:record-matching}, this proves the proposition.
\end{proof}

\section{A degree bound for matchings}\label{sec:matching}

In this section we prove Proposition~\ref{prop:matching-degree}. For a
graph $F$, let
\[
 \Psi(F)=\frac12\sum_{z\in V(F)}\sqrt{d_F(z)},
\]
so that Proposition~\ref{prop:matching-degree} states that
$M(F)\le2^{\Psi(F)}$. Equality holds when $F$ is a perfect matching.

\subsection{Proof of Proposition~\ref{prop:matching-degree}}
\label{sec:matching-proof}

The proof is by induction on $e(F)$. 
We first record how $\Psi$ changes when edges are deleted. Deleting one edge at a vertex of degree $k\ge1$ decreases the contribution of that vertex to $\Psi$ by $(\sqrt k-\sqrt{k-1})/2$. This quantity is decreasing in $k$, so it is at least $(\sqrt K-\sqrt{K-1})/2$ whenever $k\le K$. Moreover, when several edges at the same vertex are deleted one by one, each deletion decreases its contribution by at least as much as the previous one.

For real $d,s\ge1$, let
\[
 \begin{aligned}
 \alpha(d,s)&=\frac12\Bigl(\sqrt d-\sqrt{d-1}+\sqrt s-\sqrt{s-1}\Bigr),\\
 \beta(d,s)&=\frac12\Bigl(\sqrt d+\sqrt s
   +(d-1)\bigl(\sqrt s-\sqrt{s-1}\bigr)
   +(s-1)\bigl(\sqrt d-\sqrt{d-1}\bigr)\Bigr).
 \end{aligned}
\]

\begin{lemma}\label{lem:matching-edge}
For all integers $1\le s\le d$ with $(d,s)\notin\{(2,2),(3,2)\}$,
\[
 2^{-\alpha(d,s)}+2^{-\beta(d,s)}\le1.
\]
\end{lemma}

We prove Lemma~\ref{lem:matching-edge} in
Section~\ref{sec:matching-scalar}.  

\begin{proof}[Proof of Proposition~\ref{prop:matching-degree} assuming Lemma~\ref{lem:matching-edge}]
If $e(F)=0$, then $M(F)=1=2^{\Psi(F)}$. Let $e(F)\ge1$ and assume
$M(F')\le2^{\Psi(F')}$ whenever $e(F')<e(F)$. Isolated vertices affect
neither $M(F)$ nor $\Psi(F)$, and both $\log_2M$ and $\Psi$ are additive
over components, so we may assume that $F$ is connected.

Let $u$ be a vertex of maximum degree $d=\Delta(F)$, let $v$ be a
neighbor of $u$ of maximum degree among the neighbors of $u$, and let
$s=d_F(v)$. Then $1\le s\le d$, and every neighbor of $u$ has degree at
most $s$. The matchings of $F$ not containing $uv$ are exactly the
matchings of $F-uv$, and those containing $uv$ correspond to the
matchings of $F-u-v$. Hence
\begin{equation}\label{eq:matching-edge-recurrence}
 M(F)=M(F-uv)+M(F-u-v).
\end{equation}
It remains to compare the decrease of $\Psi$ in the two branches.
Deleting $uv$ decreases $\Psi$ by exactly $\alpha(d,s)$. Deleting $u$
and $v$ removes their contributions $(\sqrt d+\sqrt s)/2$ and also
deletes the remaining $d-1$ edges at $u$ and $s-1$ edges at $v$. The
former edges end at neighbors of $u$, which have degree at most $s$, and the
latter end at vertices of degree at most $d$. At a common neighbor of
$u$ and $v$ we delete the edge to $u$ first. Hence, by the monotonicity
noted above,
\begin{equation}\label{eq:matching-potential-losses}
 \Psi(F)-\Psi(F-uv)=\alpha(d,s),\qquad
 \Psi(F)-\Psi(F-u-v)\ge\beta(d,s).
\end{equation}
Since $F-uv$ and $F-u-v$ have fewer edges than $F$, the induction
hypothesis,~\eqref{eq:matching-edge-recurrence} and
\eqref{eq:matching-potential-losses} give
\[
 M(F)\le2^{\Psi(F-uv)}+2^{\Psi(F-u-v)}
 \le2^{\Psi(F)-\alpha(d,s)}+2^{\Psi(F)-\beta(d,s)},
\]
that is,
\begin{equation}\label{eq:matching-normalized-branches}
 \frac{M(F)}{2^{\Psi(F)}}\le2^{-\alpha(d,s)}+2^{-\beta(d,s)}.
\end{equation}
If the right-hand side is at most $1$, then $M(F)\le2^{\Psi(F)}$ and the
induction step is complete. By Lemma~\ref{lem:matching-edge} this holds
unless $(d,s)\in\{(3,2),(2,2)\}$, which we treat separately in the same
way, using sharper estimates.

In each case below, we apply the induction hypothesis to each smaller
graph, bound its matching count by $2^{\Psi}$, and use the corresponding
lower bound on the decrease of $\Psi$ before dividing by $2^{\Psi(F)}$.

\medskip\noindent\textbf{Case 1: $(d,s)=(3,2)$.}
Let $w$ be the neighbor of $v$ other than $u$. Every neighbor of $u$ has
degree at most $2$.

\smallskip\noindent\emph{Case 1(a): $d_F(w)\le2$.}
Then the three edges from $\{u,v\}$ to the rest of $F$ all end at
vertices of degree at most $2$, and each of their deletions decreases
$\Psi$ by at least $(\sqrt2-1)/2$. Hence
\[
 \Psi(F)-\Psi(F-u-v)\ge\frac{\sqrt3+\sqrt2+3(\sqrt2-1)}2
 =\frac{\sqrt3+4\sqrt2-3}2,
\]
and since $\alpha(3,2)=(\sqrt3-1)/2$, the induction hypothesis gives
\[
 \frac{M(F)}{2^{\Psi(F)}}
 \le2^{-\alpha(3,2)}+2^{-(\Psi(F)-\Psi(F-u-v))}
 \le2^{-(\sqrt3-1)/2}+2^{-(\sqrt3+4\sqrt2-3)/2}<0.995<1.
\]

\smallskip\noindent\emph{Case 1(b): $d_F(w)=3$.}
Then $w$ is not adjacent to $u$. Classifying matchings according to
whether $v$ is unmatched, matched to $u$ or matched to $w$,
\begin{equation}\label{eq:matching-vertex-recurrence}
 M(F)=M(F-v)+M(F-v-u)+M(F-v-w).
\end{equation}
The induction hypothesis applies to $F-v$, $F-v-u$ and $F-v-w$.
Deleting $v$ decreases $\Psi$ by
$\frac{\sqrt2}2+2\cdot\frac{\sqrt3-\sqrt2}2=\sqrt3-\frac{\sqrt2}2>1$,
so $M(F-v)\le2^{\Psi(F)-1}$.
In $F-v$, the vertices $u$ and $w$ have degree $2$ and every vertex has
degree at most $3$. Hence, for $x\in\{u,w\}$, deleting $x$ from $F-v$
decreases $\Psi$ by at least
$\frac{\sqrt2}2+2\cdot\frac{\sqrt3-\sqrt2}2=\sqrt3-\frac{\sqrt2}2>1$.
Thus, in total $\Psi(F)-\Psi(F-v-x)>2$ for $x\in\{u,w\}$, hence
$M(F-v-x)\le2^{\Psi(F)-2}$.
Using~\eqref{eq:matching-vertex-recurrence} and dividing by
$2^{\Psi(F)}$,
\[
 \frac{M(F)}{2^{\Psi(F)}}<\frac12+\frac14+\frac14=1.
\]

\medskip\noindent\textbf{Case 2: $(d,s)=(2,2)$.}
Since $F$ is connected and $\Delta(F)=2$, it is a path or a cycle.

\smallskip\noindent\emph{Case 2(a): $F$ is a path.}
Then $F$ has at least three vertices. Let $xy$ be an end edge with
$d_F(x)=2$ and $d_F(y)=1$. Deleting $xy$ decreases $\Psi$ by
$\frac{\sqrt2-1}2+\frac12=\frac1{\sqrt2}$, and deleting $x$ and $y$
decreases it by at least
$\frac{\sqrt2+1}2+\frac{\sqrt2-1}2=\sqrt2$. The induction hypothesis applies to $F-xy$ and $F-x-y$. Using the
recurrence at $xy$ and the decrease bounds above,
\[
 \frac{M(F)}{2^{\Psi(F)}}\le2^{-1/\sqrt2}+2^{-\sqrt2}<0.988<1.
\]

\smallskip\noindent\emph{Case 2(b): $F$ is a cycle.}
Let $F=C_m$, so that $\Psi(F)=m/\sqrt2$. For $m=3,4$ we have $M(C_3)=4$
and $M(C_4)=7$, and
\[
 4\cdot2^{-3/\sqrt2}<0.920,\qquad 7\cdot2^{-2\sqrt2}<0.986.
\]
For $k\ge3$, the recurrence at
an edge of $C_k$ gives $M(C_k)=M(P_k)+M(P_{k-2})$, where   $P_k$ is the path with $k$ vertices.  The recurrence at
an end edge of $P_k$ gives $M(P_k)=M(P_{k-1})+M(P_{k-2})$. Together they
give $M(C_k)=M(C_{k-1})+M(C_{k-2})$ for $k\ge5$. For $m\ge5$, the
recurrence and induction hypothesis applied to the smaller graphs
$C_{m-1}$ and $C_{m-2}$ give
\[
 \frac{M(C_m)}{2^{m/\sqrt2}}\le2^{-1/\sqrt2}+2^{-\sqrt2}<0.988<1.
\]
This completes the induction.
\end{proof}

\subsection{Proof of Lemma~\ref{lem:matching-edge}}
\label{sec:matching-scalar}

For $(d,s)=(1,1)$ we have $\alpha=\beta=1$, and the inequality holds
with equality. We prove the strict inequality for all other pairs. For
real $d,s\ge1$, let
\[
 R(d,s)=\frac{2^{-\beta(d,s)}}{1-2^{-\alpha(d,s)}}.
\]
Since $\alpha(d,s)>0$, the strict inequality is equivalent to
$R(d,s)<1$. The functions $\alpha,\beta$ and $R$ are symmetric in $d$
and $s$. We show that $R(d,s)$ is decreasing in $d$ for
$d\ge\max\{2,s\}$; by symmetry, $R(x,x)$ is then decreasing for $x\ge2$.
This reduces the lemma to the pairs $(2,1)$, $(4,2)$ and $(3,3)$.

Fix $s\ge1$ and let $d\ge\max\{2,s\}$. Writing
$\eta=\sqrt d-\sqrt{d-1}$, a direct computation gives
\[
 \partial_d\alpha(d,s)=-\frac{\eta}{4\sqrt{d(d-1)}},\qquad
 \partial_d\beta(d,s)
 =\frac{\sqrt s-\sqrt{s-1}}2
  +\frac{\eta}4\Bigl(1+\frac{d-s}{\sqrt{d(d-1)}}\Bigr)
 \ge\frac34\eta,
\]
where the inequality uses $d\ge s$ and
$\sqrt s-\sqrt{s-1}\ge\eta$. Moreover, $\alpha(d,s)\ge\eta$, and
$2^x-1\ge x\ln2$ for $x>0$. Since
$1/\sqrt{d(d-1)}=\eta\bigl(1/\sqrt d+1/\sqrt{d-1}\bigr)$, we obtain
\[
 -\partial_d\log_2R(d,s)
 =\partial_d\beta(d,s)+\frac{\partial_d\alpha(d,s)}{2^{\alpha(d,s)}-1}
 \ge\frac{\eta}4\left(3-\frac1{\ln2}
   \Bigl(\frac1{\sqrt d}+\frac1{\sqrt{d-1}}\Bigr)\right)>0,
\]
where the last inequality uses $d\ge2$, $1/\ln2<3/2$ and
$1/\sqrt d+1/\sqrt{d-1}\le1+1/\sqrt2<2$.
Thus $R(d,s)$ is decreasing in $d$ for $d\ge\max\{2,s\}$, and
\[
 \frac{d}{dx}\log_2R(x,x)
 =2\,\partial_d\log_2R(d,s)\big|_{d=s=x}<0\qquad(x\ge2).
\]

Direct computation gives
\[
 2^{-\alpha(2,1)}+2^{-\beta(2,1)}<0.919,\qquad
 2^{-\alpha(4,2)}+2^{-\beta(4,2)}<0.971,\qquad
 2^{-\alpha(3,3)}+2^{-\beta(3,3)}<0.997.
\]
For $s=1$, monotonicity in $d$ gives the result for all $d\ge2$. For
$s=2$, it gives the result for all $d\ge4$. For $3\le s\le d$,
\[
 R(d,s)\le R(s,s)\le R(3,3)<1.
\]
This proves Lemma~\ref{lem:matching-edge}.\qed

\section*{Declaration of AI use}  
The author provided the main ideas of the proof, combining the hypergraph container method with the deferral of diagonal constraints. GPT-6 Pro assisted in completing the proofs, and GPT 6.1 and Opus 5.5 were used to revise the manuscript and improve its exposition. The author has verified all mathematical arguments and results and takes full responsibility for the content of this paper.

\end{document}